\documentclass[11pt]{article}
\usepackage[english]{babel}
\usepackage[latin1]{inputenc}
\usepackage{exscale}
\usepackage[centertags]{amsmath}
\usepackage{amsfonts,amstext,amssymb,amsthm}
\usepackage{comment}
\usepackage{newlfont}
\usepackage{enumerate}
\usepackage{graphicx}
\usepackage{HolderContinuityAbraham2}
\usepackage{color}
\usepackage[usenames,dvipsnames]{xcolor}
\definecolor{BlueFonse}{rgb}{0,0,1}
\definecolor{BlueFonse1}{cmyk}{1,0,0,0.7}
\usepackage[colorlinks=true,linkcolor=blue,citecolor=ForestGreen]{hyperref}
\usepackage{hyperref}
\usepackage{url}
\title{Local H\"older continuity of the isoperimetric profile in complete noncompact Riemannian manifolds with bounded geometry}
\author{Abraham Mu\~noz Flores\footnote{Partially supported by Capes}, Stefano Nardulli}
\begin{document}
      \maketitle
\pretolerance=2000
\tolerance=3000
\noindent {\sc abstract}. For a complete noncompact connected Riemannian manifold with bounded geometry $M^n$, we prove that the isoperimetric profile function $I_{M^n}$ is a locally $\left(1-\frac{1}{n}\right)$-H\"older continuous function and so in particular it is continuous. Here for bounded geometry we mean that $M$ have $Ricci$ curvature bounded below and volume of balls of radius $1$, uniformly bounded below with respect to its centers. We prove also the equivalence of the weak and strong formulation of the isoperimetric profile function in complete Riemannian manifolds which is based on a lemma having its own interest about the approximation of finite perimeter sets with finite volume by open bounded with smooth boundary ones of the same volume. Finally the upper semicontinuity of the isoperimetric profile for every metric (not necessarily complete) is shown.
\bigskip\bigskip

\bigskip

\centerline{\bf AMS subject classification: }
49Q20, 58E99, 53A10, 49Q05.
      \tableofcontents       
\section{Introduction}\label{1} 
In this paper we always assume that all the Riemannian manifolds $(M, g)$ considered are smooths with smooth Riemannian metric $g$. We denote by $V_g$ the canonical Riemannian measure induced on $M$ by $g$, and by $A_g$ the $(n-1)$-Hausdorff measure associated to the canonical Riemannian length space metric $d$ of $M$, by $\P_g(\Omega, U)$ the perimeter in $U\subseteq M$ with respect to the metric $g$ of a finite perimeter set $\Omega\subseteq M$, here $U$ is an open set, by $|Du|_g$ we denote the positive Radon measure represented by the total variation of the distributional gradient of a $BV$-function $u$ having domain $M$.  For each $k\in\R$ we denote by $\mathbb{M}_k^n$ the $n$-dimensional space form of constant sectional curvature equal to $k$. When it is already clear from the context, explicit mention of the metric $g$ will be suppressed. When dealing with finite perimeter sets or locally finite perimeter sets we will denote the reduced boundary by $\partial^*\Omega$, whenever the topological boundary $\partial\Omega$ is smooth the reduced boundary coincides with the topological boundary $\partial\Omega$. For this reason we will denote $\P(\Omega):=\P(\Omega, M)=A(\partial^*\Omega)=A(\partial\Omega)$ when no confusion may rise, and for every finite perimeter set $\Omega'$ we always choose a representative $\Omega$ (i.e., that differs from $\Omega'$ by a set of Riemannian measure $0$), such that $\partial\Omega=\overline{\partial^*\Omega}$, where $\partial\Omega$ is the topological boundary of $\Omega$. At this point we give the definition of the isoperimetric profile function which is our main object of study in this paper.   
\subsection{The isoperimetric profile}
\begin{Def}\label{Def:IsPStrong}
Typically in the literature, the \textbf{isoperimetric profile function of $M$} (or briefly, the isoperimetric profile) $I_M:[0,V(M)[\rightarrow [0,+\infty [$, is defined by $$I_M(v):= \inf\{A(\partial \Omega): \Omega\in \tau_M, V(\Omega )=v \},$$ where $\tau_M$ denotes the set of relatively compact open subsets of $M$ with smooth boundary.
\end{Def}
However there is a more general context in which to consider this notion that will be better suited to our purposes.  Namely, we can give a weak formulation of the preceding variational problem replacing the set $\tau_M$ with the family $\tilde{\tau}_M$ of subsets of finite perimeter of $M$. 
\begin{Def}\label{Def:IsPWeak}
Let $M$ be a Riemannian manifold of dimension $n$ (possibly with infinite volume). We denote by $\tilde{\tau}_M$ the set of  finite perimeter subsets of $M$. The function $\tilde{I}_M:[0,V(M)[\rightarrow [0,+\infty [$  defined by 
     $$\tilde{I}_M(v):= \inf\{\mathcal{P}(\Omega): \Omega\in \tilde{\tau}_M, V(\Omega )=v\},$$ 
is called the \textbf{weak isoperimetric profile function} (or shortly the \textbf{isoperimetric profile}) of the manifold $M$. If there exists a finite perimeter set $\Omega\in\tilde{\tau}_M$ satisfying $V(\Omega)=v$, $\tilde{I}_M(V(\Omega))=A(\partial^*\Omega)= \mathcal{P}(\Omega)$ such an $\Omega$ will be called an \textbf{isoperimetric region}, and we say that $\tilde{I}_M(v)$ is \textbf{achieved}. 
\end{Def} 
There are many others possible definitions of isoperimetric profile corresponding to the minimization over various differents sets of admissible domains, as stated in the following definition.  
\begin{Def} For every $v\in[0, +\infty[$, let us define
\begin{eqnarray*}
I^{*}_{M}(v):=inf\{A(\partial_{top}\Omega): \Omega\subset M, \partial_{top}\Omega\: \text{is}\: C^{\infty}, V(\Omega)=v\},\\
\tilde{I}_{M}^{*}(v):=inf\{\mathcal{P}_{M}(\Omega):\Omega\subset M, \Omega\in\tilde{\tau}_{M} , V(\Omega)=v, diam(\Omega)<+\infty\}, 
\end{eqnarray*}
where $diam(\Omega):=\sup\{d(x,y):x,y\in\Omega\}$ denotes the \textbf{diameter} of $\Omega$.
\end{Def}
\begin{Rem}\label{Rem:Trivialinequality}
Trivially one have $I_M\geq I^{*}_{M}\geq\tilde{I}_M$ and $I_M\geq\tilde{I}^*_{M}\geq\tilde{I}_M$. 
\end{Rem}
However as we will see in Theorem $\ref{Thm:Equivalence}$, all of these definitions are actually equivalents, in the sense that the infimum remains unchanged, i.e., $I_M=\tilde{I}_M$.
\subsection{Main Results}\label{MainRes} 
\begin{Res}\label{Thm:Equivalence} If $M^n$ is an arbitrary complete Riemannian manifold, then $I_M(v)=\tilde{I}_{M}^{*}(v)=\tilde{I}_M(v)=I^{*}_{M}(v)$. 
\end{Res} 
The proof of this fact involves actually very natural ideas. In spite of this it is technical and we have found no written traces in the literature, unless Lemma $2$ of \cite{Modica} that deal with the case of a compact domain of $\R^n$ as an ambient space. Hence we provided ourselves a proof based on Lemma \ref{Lemma:Smoothisovolumic} which have an independent interest, because it gives an approximation theorem of a finite perimeter set by open relatively compact sets with smooth boundary of the same volume and for this constitutes a refinement of a more classical approximation theorem of finite perimeter sets by members of $\tau_M$ that one can find in the literature (see for example the books of \cite{Maggi}, \cite{AFP},  or in the paper \cite{MPPP}).  The equivalence stated in Theorem \ref{Thm:Equivalence} allows us to consider elements of $\tau_M$ or $\tilde{\tau}_M$ according to what is more convenient in subsequent arguments. This observation is used in a crucial way when we prove Theorem $\ref{Main}$ and Corollary \ref{CorRes:Uppersemicontinuity}. This latter could be considered as a corollary of Lemma \ref{Lemma:HoleAproximation}.
\begin{CorRes}\label{CorRes:Uppersemicontinuity} Let $M^n$ be an $n$-dimensional Riemannian manifold (possibly incomplete, or possibly complete not necessarily with bounded geometry). Then $I_M$ is upper semicontinuous.  
\end{CorRes}
\begin{Def}\label{Def:BoundedGeometry}
A complete Riemannian manifold $(M, g)$, is said to have \textbf{bounded geometry} if there exists a constant $k\in\mathbb{R}$, such that $Ric_M\geq k(n-1)$ (i.e., $Ric_M\geq (n-1)kg$ in the sense of quadratic forms) and $V(B_{(M,g)}(p,1))\geq v_0$ for some positive constant $v_0$, where $B_{(M,g)}(p,r)$ is the geodesic ball (or equivalently the metric ball) of $M$ centered at $p$ and of radius $r> 0$.
\end{Def}
\begin{Res}[Local $\left(1-\frac{1}{n}\right)$-H\"older continuity of the isoperimetric profile]\label{MainLemma}
Let $M^n$ be a complete smooth Riemannian manifold with bounded geometry. Then there exists a positive constant $C=C(n,k)$ such that for every $v,v'\in]0, V(M)[$ satisfying $|v-v'|\leq \frac{1}{C(n,k)}\min\left(v_0,\left(\frac{v}{I_M(v)+C(n,k)}\right)^n\right)$, we have  
\begin{equation}\label{Eq:ResMainLemmaStatement}
\left|I_M(v)-I_M(v')\right|\leq C(n,k)\left(\frac{|v-v'|}{v_0}\right)^{\frac{n-1}{n}}.
\end{equation}
In particular $I_M$ is continuous on $[0,V(M)[$.
\end{Res}
\begin{Def} Let $f:(X,d)\to\R$ and $\alpha\in[0,1]$, we say that $f$ is \textbf{locally $\alpha$-H\"older continuous on $X$}, for every $z\in X$ there exist $\delta_z, C_z>0$ such that for every $x,y\in X$ satisfying $|x-z|, |y-z|\leq\delta_z$ we have $|f(x)-f(y)|\leq C_z|x-y|^{\alpha}$.We say that $f$ is \textbf{uniformly locally $\alpha$-H\"older continuous on $X$}, if there exist two constants $\delta, C>0$ such that for every $x,y\in X$ satisfying $d(x,y)\leq\delta$ we have $|f(x)-f(y)|\leq C|x-y|^{\alpha}$. We say that $f$ is \textbf{(globally) $\alpha$-H\"older continuous on $X$}, if there exists $C>0$ such that $|f(x)-f(y)|\leq C|x-y|^{\alpha}$ for every $x,y\in X$. We call the various constants $C_z, C$ appearing in this definition the \textbf{H\"older constants of $f$}.
\end{Def}
\begin{CorRes}[Local $\frac{n-1}{n}$-H\"older continuity of the isoperimetric profile]\label{Main}
Let $M^n$ be a complete smooth Riemannian manifold with bounded geometry and $v\in]0, V(M)[$. Then there exists positive constants $\delta=\delta(n,k,v_0,v)>0$, if $k\leq 0$, $\delta=\delta(n,k,v_0,, v, V(M))$, if $k>0$, and $C=C(n,k)>0$, such that for every $v_1,v_2\in[v-\delta, v+\delta]$ we have 
\begin{equation}\label{Eq:ResMainStatement}
\left|I_M(v_1)-I_M(v_2)\right|\leq C(n,k)\left(\frac{|v_1-v_2|}{v_0}\right)^{\frac{n-1}{n}}.
\end{equation}
Moreover, if $V(M)=+\infty$ then $I_M$ is uniformly locally $\frac{n-1}{n}$-H\"older continuous on $[\bar{v}, +\infty[$, for every $\bar{v}>0$. If $V(M)=+\infty$ then $I_M$ is globally $\frac{n-1}{n}$-H\"older continuous on every interval $[a,b]\subset ]0, +\infty[$ with H\"older constant $\bar{C}$ depending on $n,k,v_0, a,b$. If $V(M)<+\infty$, then $I_M$ is globally $\frac{n-1}{n}$-H\"older continuous on $[\bar{v},V(M)-\bar{v}]$, for every $\bar{v}\in]0,\frac{V(M)}{2}[$.
\end{CorRes}
\begin{Rem} Unfortunately $\lim_{a\to 0^+}\bar{C}(n,k,v, a,b)=+\infty$ and 

$\lim_{b\to 0^+}\bar{C}(n,k,v, a,b)=+\infty$.
\end{Rem}
\begin{Rem} Observe that in the statement of the preceding Corollary the H\"older constant $C$ does not depend on $v_0$ and $v$, but just $\delta$ depends on them.
\end{Rem}
\begin{Rem} At our actual knowledge, it is still an open question wether or not we can prove global $\frac{n-1}{n}$-H\"older continuity of $I_M$ on an arbitrary proper interval $[0,b]\subset[0, V(M)[$ or on the entire interval $[0,V(M)[$, or at least unifom local $\frac{n-1}{n}$-H\"older continuity on $[0, V(M)[$, when we assume the manifold $M$ to be with bounded geometry and with $V(M)=+\infty$. 
\end{Rem}
 The next fact to be observed is that it is worth to have a proof of the continuity or H\"older continuity of the isoperimetric profile, because in general the isoperimetric profile function of a complete Riemannian manifold is not continuous. In case of manifolds with density, in Proposition $2$ of \cite{MorganBlog}  is exhibited an example of  a manifold with density having discontinuous isoperimetric profile. To exhibit a complete Riemannian manifold with a discontinuous isoperimetric profile is a more  subtle and difficult task that was performed by the second author and Pierre Pansu in \cite{NardulliPansuDiscontinuous}, for manifolds of dimension $n\geq 3$.  In spite of these quite sophisticated counterexamples the class of manifolds admitting a continuous isoperimetric profile is vast, for an account of the existing literature on the continuity results obtained for $I_M$, one could consult the introduction of \cite{RitoreContinuity} and the references therein. If $M$ is compact, classical compactness arguments of geometric measure theory combined with the direct method of the calculus of variations provide a short proof of the continuity of $I_{M}$ in any dimension $n$, \cite{MorganBlog} Proposition $1$. Finally, if $M$ is complete, non-compact, and $V(M)<+\infty$, an easy consequence of Theorem $2.1$ in \cite{RRosales} yields the possibility of extending the same compactness argument valid in the compact case and to prove the continuity of the isoperimetric profile, see for instance Corollary 2.4 of \cite{NardulliRusso}. A careful analysis of Theorem $1$ of \cite{Nar12} about the existence of generalized isoperimetric regions, leads to the continuity of the isoperimetric profile $I_M$ in manifolds with bounded geometry satisfying some other assumptions on the geometry of the manifold at infinity, of the kind considered by the second author and A. Mondino in \cite{MonNar}, i.e., for every sequence of points diverging to infinity, there exists a pointed smooth manifold $(M_{\infty}, g_{\infty}, p_{\infty})$ such that $(M,g,p_j)\rightarrow (M_{\infty}, g_{\infty}, p_{\infty})$ in $C^0$-topology. This proof is independent from that of Theorem \ref{Main}. This is not the case for general complete infinite-volume manifolds $M$. Recently Manuel Ritor\'e (see for instance \cite{RitoreContinuity}) showed that a complete Riemannian manifold possessing a strictly convex Lipschitz continuous exhaustion function has continuous  and nondecreasing isoperimetric profile $\tilde{I}_M$. Particular cases of these manifolds are Cartan-Hadamard manifolds and complete noncompact manifolds with strictly positive sectional curvatures. In \cite{RitoreContinuity} as in our Theorem \ref{Main} the major difficulty consists in finding a suitable way of subtracting a volume to an almost minimizing region. 

The aim of this paper is to prove Theorem \ref{Main} in which we give a very short and quite elementary proof of the continuity of $I_M$ when $M$ is a complete noncompact Riemannian manifold of bounded geometry and even better we show that $I_M$ is actually a locally $C^{1-\frac{1}{n}}(]0,V(M)[)$ function. The reason which allow us to achieve this goal, is that in bounded geometry it is always possible to add or subtract to a finite perimeter set a small ball that captures a fixed fraction of volume (depending only from the bounds of the geometry) centered at points close to it. Corollary \ref{CorRes:Uppersemicontinuity} ensures upper semicontinuity, so the problems appears when we try to prove lower semicontinuity. To prove lower semicontinuity we need some kind of compactness that is expressed here by a bounded geometry condition.  Geometrically speaking our assumptions of bounded geometry ensures that the manifold at infinity is thick enough to permit to place a small geodesic ball $B$ close to an arbitrary domain $D$ in such a way $V(B\cap D)$ recovers a controlled fraction of $V(D)$ and this fraction depends only on $V(D)$ and the bounds on the geometry $n, v_0, k$, see Definition \ref{Def:BoundedGeometry} for the exact meaning of $n$, $v_0$, $k$.  The proof that we present here uses only metric properties of the manifolds with bounded geometry and for this reason it is still valid when suitably reformulated in the context of metric measure spaces. One can find similar ideas already in the metric proof of continuity of the isoperimetric profile contained in \cite{Gallot}. For the full generality of the results we need that the spaces have to be doubling, satisfying a $1$-Poincar\'e inequality and a curvature dimension condition. This class of metric spaces includes for example manifolds with density as well  as subRiemannian manifolds. Following the arguments contained in \cite{BP} we can obtain another proof of the continuity of the isoperimetric profile under our assumptions of bounded geometry but with the extra assumption of the existence of isoperimetric regions of every volume, which is less general of our own proof of Theorem $\ref{Main}$, because in Theorem $\ref{Main}$ we do not need to assume any kind of existence of isoperimetric regions. In spite of this the Heintze-Karcher type arguments used in \cite{BP} have an advantage because they permits to give a uniform bound on the length of the mean curvature vector of the generalized isoperimetric regions (i.e., left and right derivatives of $I_M$) with volumes inside an interval $[a, b]\subset ]0, V(M)[$, depending only on $a$ and $b$. Finally, we mention that just with Ricci bounded below and existence of isoperimetric regions the arguments of \cite{BP} fails and we cannot prove the continuity of the isoperimetric profile, for this we need a noncollapsing condition on the volume of geodesic balls as in our definition of bounded geometry. 
\begin{Rem} It remains still an open question whether $Ricci$ bounded below and existence of isoperimetric regions for every volume implies continuity of the isoperimetric profile in presence of collapsing. We are not able to extend to this setting the arguments of $\cite{BP}$, neither to provide a counterexample, because the manifolds with discontinuous isoperimetric profile constructed in $\cite{NardulliPansuDiscontinuous}$ have $Ricci$ curvature tending to $-\infty$.
\end{Rem}
\subsection{Plan of the article}
\begin{enumerate}
           \item  Section $\ref{1}$ constitutes the introduction of the paper. We state the main results of the paper.
           \item In Section $\ref{Sec:Weak}$ we prove that $\tilde{I}_M=I_M$.
           \item In section $\ref{Sec:Continuity}$ we prove the local $C^{1-\frac{1}{n}}$-H\"older continuity of the isoperimetric profile in bounded geometry, i.e., Theorem $\ref{MainLemma}$ and Corollary $\ref{Main}$ without assuming existence of isoperimetric regions. 
\end{enumerate}
\subsection{Acknowledgements}  
The second author is indebted to Pierre Pansu for inspiring this paper and then to Pierre Pansu, Frank Morgan, Andrea Mondino, and Luigi Ambrosio for useful discussions on the topics of this article. 
The first author wish to thank the CAPES for financial support for the period in which he was a Ph.D. student at IM-UFRJ. Finally we want to thank a lot the anonymous referee whose comments contributed to improve both the results and the presentation of the proofs contained in this paper.   

 \section{Equivalence of the weak and strong formulation}\label{Sec:Weak} 
 As the example $3.53$ of \cite{AFP} shows, in general we can have finite perimeter sets with positive perimeter and void interior that are not equivalent to any other set of finite perimeter with non void interior. So the question of putting a ball inside or outside a set of finite perimeter is a genuine technical problem. On the other hand, following \cite{Tamanini} Theorem $1$, it is always possible to put a small ball inside and outside an isoperimetric region. As a general remark a result of Federer (the reader could consult \cite{AFP} Theorem $3.61$) states that for a given set of finite perimeter $E$ the density is either $0$ or $\frac{1}{2}$ or $1$, $\mathcal{H}^{n-1}$-a.e. $x\in M$, moreover points of density $1$ always exist $V$-a.e. inside $D$, because of the Lebesgue's points Theorem applied to the characteristic function of any $V$-measurable set of $M$. About this topic the reader could consult the book \cite{Maggi} Example $5.17$. Thus $V(D)>0$ ensures the existence of at least one point $p$ belonging to $D$ of density $1$, which is enough for the aims of our proofs. In view of these facts to prove Theorem $\ref{Thm:Equivalence}$ we need to make a construction which replace a finite perimeter set by one of the same volume with a small ball inside and one outside, by adding a small geodesic ball (with smooth boundary) to a point of density $0$ and subtracting a small geodesic ball to a point of density $1$ taking care of not altering the volume. This enables us to obtain again a finite perimeter set of the same volume with a perimeter that is a small perturbation of the original one and that in addition have the property that we can put inside and outside a small ball. This construction legitimate us to apply mutatis mutandis the arguments of the proof of Lemma $1$ of \cite{Modica} to get the isovolumic approximation Lemma \ref{Lemma:Smoothisovolumic} and then to conclude the proof of Theorem $\ref{Thm:Equivalence}$. Our adapted version of Lemma $1$ of \cite{Modica} is the following lemma. 
\begin{Lemme}\label{Lemma:1ModicaReferee} Let $\Omega_1\in\tilde{\tau}_M$ with $V(\Omega_1)<+\infty$, such that there exists two geodesic balls satisfying $B(x_1, r_1)\subset\Omega_1$ and $B(x_2,r_2)\cap\Omega_1=\emptyset$, with $0<r_1<inj_M(x_1)$ and $0<r_2<inj_M(x_2)$. We set $v^*:=\min\left\{V(B(x_1, \frac{r_1}{2})), V(B(x_2, \frac{r_2}{2}))\right\}$. For any $v\in[0,v^*]$ we denote by $R_{i,v}$ a radius such that $V(B(x_i,R_{i,v}))=v$ and by $S(x,r)$ the sphere of radius $r$ and center $x$. Let us define 
\begin{equation}\label{Eq:defoff}
f_{\Omega_1}(v):=\max\left\{\sup_{0\leq t\leq R_{1,v}}A(S(x_1,t)), \sup_{0\leq t\leq R_{2,v}}A(S(x_2,t))\right\}.
\end{equation} 
Then for any $\varepsilon>0$ and any $v\in]V(\Omega_1)-v^*, V(\Omega_1)+v^*[$, there exists $\Omega_2\in\tau_M$ such that $V(\Omega_2)=v$ and
 $$\P(\Omega_2)\leq\P(\Omega_1)+f_{\Omega_1}(|v-V(\Omega_1)|)+\frac{\varepsilon}{4}.$$ 
\end{Lemme}
\begin{Rem}
We observe that if $M$ is noncompact and $\Omega$ bounded, then we always have $Interior(\Omega^c)\neq\emptyset$. 
\end{Rem}
\begin{Dem}[of Lemma \ref{Lemma:1ModicaReferee}]
By the proof of the claim p. 105 of \cite{MPPP}, there exists a sequence of $BV$-functions $(u_l)$ on $M$ such that $\lim_l||u_l-\chi_{\Omega_1}||_1=0$, $|Du_l|(M)=\P(\Omega_1)$ and each $u_l$ has compact support $K_l$. Note that we can assume that $B(x_1, r_1)\subset K_l$. Moreover, construction the $u_l$ satisfy $0\leq u_l\leq\chi_{\Omega_1}$, which gives $K_l\subset\Omega_1$. Considering a smooth positive kernel $\rho$ with compact support the mollified functions $u_{j,l}=u_l*\rho_{\frac{1}{j}}$ satisfy $0\leq u_{j,l}\leq 1$, $\lim_{j\to+\infty}||u_{j,l}-u_l||_1=0$, $\lim_{l}|Du_{j,l}|(M)=|Du_l|(M)$ and for $j$ large enough the support $K_{j,l}$ of $u_{j,l}$ satisfies $B(x_1,\frac{r}{2})\cap K_{j,l}=\emptyset$.
\begin{Rem} As explained in \cite{MPPP} to perform a convolution on a manifold one have just to use a partition of unity associated to finite sets of local charts covering the compact support of $u_l$ and then mollify in each local chart. 
\end{Rem}
By a diagonal argument we extract a subsequence $v_l=u_{j,l}$, satisfying $0\leq v_l\leq 1$, $\lim_{l}||v_l-\chi_{\Omega_1}||_1=0$, $\lim_l|Dv_l|(M)=\P(\Omega_1)$, and for $l$ large enough the support $C_l$ of $v_l$ satisfies $B(x_1,\frac{r_1}{2})\subset C_l$ and $B(x_2,\frac{r_2}{2})\cap C_l=\emptyset$. Putting $F_t^l:=\{x\in M:v_l(x)>t\}$ and using the Fleming-Rishel Theorem (compare Theorem $4.3$ of \cite{ADo}) we have 
\begin{equation*}
\P(\Omega_1)=\lim_l|Dv_l|=\lim_l\int_0^1\P(F_t^l)dt\geq\int_0^1\liminf_l\P(F_t^l)dt.
\end{equation*}
An application of Sard's Theorem ensures that the sets $F_t^l$ are smooth for almost every $t\in]0,1[$. Thus for every $l$ we can choose a $t\in]0,1[$ (depending on $l$), such that $\liminf_l\P(F_t^l)\leq\P(\Omega_1)$. Moreover, we have $|V(F_t^l)-V(\Omega_1)|\leq V(F_t^l\setminus\Omega_1)+V(\Omega_1\setminus F_t^l)$ and 
\begin{equation*}
V(F_t^l\setminus\Omega_1)\leq\frac{1}{t}||v_l-\chi_{\Omega_1}||_1,
\end{equation*} 
\begin{equation*}
V(\Omega_1\setminus F_t^l)\leq\frac{1}{1-t}||v_l-\chi_{\Omega_1}||_1.
\end{equation*}
Since we have $|v-V(\Omega_1)|<v_0$, we can choose $l$ large enough to get 
\begin{equation*}
|v-V(\Omega_1)|+\frac{||v_l-\chi_{\Omega_1}||_1}{t(1-t)}<v^*,
\end{equation*}
which yields for $l$ large enough $|V(F_t^l)-v|<v^*$. Hence by subtracting $B(x_1, R_{1,V(F_t^l)-v})$ or adding $B(x_2, R_{2,v-V(F_t^l)})$ to $F_t^l$, we obtain a bounded open set with smooth topological boundary and volume $v$ and perimeter equal to 
\begin{equation*}
\P(F_t^l)+A(S(x_{i,l}, R_{i,l})\leq\P(F_t^l)+f_{\Omega_1}(|v-V(F_t^l)|),
\end{equation*}
where $R_{i,l}:=R_{2,v-V(F_t^l)}$ if $V(F_t^l)<v$ and $R_{i,l}:=R_{1,V(F_t^l)-v}$, if $V(F_t^l)<v$ and $R_{i,l}=0$ if $V(F_t^l)=v$ otherwise. We finally get $\Omega_2$ for any $l$ large enough and we conclude the proof.
\end{Dem}

We can state now the next lemma which permits to approximate an arbitrary finite perimeter set with another one having the same volume and two holes (balls), one inside and the other outside it. Before stating the next lemma just let us mention that for a set $X$ inside a topological space we denote by $Interior(X)=\mathring{X}$ the set of its interior points. 
\begin{Lemme}\label{Lemma:HoleAproximation} Let $M$ be a Riemannian manifold and $\Omega\in\tilde{\tau}_M$ be a set of finite perimeter with finite volume $V(\Omega)\in]0, V(M)[$. For any $\varepsilon>0$, there exists a set of finite perimeter $\tilde{\Omega}\subseteq M$ and two geodesic balls $B(x_1, r_1)$, and $B(x_2, r_2)$ such that $V(\Omega)=V(\tilde{\Omega})$, $B(x_1, r_1)\subset\Omega_1$, $B(x_2, r_2)\cap\tilde{\Omega}=\emptyset$, and 
\begin{equation}
\P(\tilde{\Omega})\leq \P(\Omega)+\frac{\varepsilon}{4}.
\end{equation}
\end{Lemme}
\begin{Dem} Consider an arbitrary set $\Omega\in\tilde{\tau}_M$ and take two distinct points $x_1\in\Omega$ and $x_2\in\Omega^c$ of density $\Theta(x_1, V\llcorner\Omega)=1$ and $\Theta(x_2, V\llcorner\Omega)=0$, where $\Theta(p, V\llcorner\Omega):=\lim_{r\to0^+}\frac{V(\Omega\cap B(p, r))}{\omega_nr^n}$, for every $p\in M$.  By $\omega_n$ we denote the volume of the ball of radius $1$ in $\R^n$. Consider the two continuous functions $f_1, f_2:I\to\R$, where $I:=[0, r_0[$ such that $f_1(r):=V(\Omega\cap B_M(x_1, r))$, $f_2(r):=V(\Omega^c\cap  B_M(x_2, r))$. The radius $r_0$ could be chosen small enough to have $B_M(x_1, r_1)\cap B_M(x_2, r_2)=\emptyset$ for every $r_1, r_2\in I$ and such that there exist $r_1, r_2\in I$ satisfying the property $f_1(r_1)=f_2(r_2)$ and $\partial B_M(x_1, r_1),\partial B_M(x_2, r_2)$ smooths (for this last property it is enough to take $r_0$ less than the injectivity radius at $x_1$ and $x_2$). Then we set 
$$\tilde{\Omega}:=[\Omega\setminus B_M(x_1,r_1)]\mathring{\cup}[\Omega^c\cap B_M(x_2, r_2)]=[\Omega\setminus B_M(x_1,r_1)]\cup B_M(x_2, r_2).$$ 
As it is easy to see $V(\tilde{\Omega})=V(\Omega)$, 
\begin{equation}\label{Eq:Equivalence}
|\P(\tilde{\Omega})-\P(\Omega)|\leq \sum_{i=1}^2[A(\partial B_M(x_i, r_i))+\P(\Omega, B_M(x_i, r_i))],
\end{equation}
\begin{equation}\label{Eq:Equivalence1}
V(\Omega\Delta\tilde{\Omega})=f_1(r_1)+f_2(r_2),
\end{equation} $\mathring{\tilde{\Omega}}\neq\emptyset$, and $Interior(\tilde{\Omega}^c)\neq\emptyset$. It is straightforward to verify that the right hand sides of \eqref{Eq:Equivalence} and \eqref{Eq:Equivalence1} converge to zero when the radii $r_1$ and $r_2$ go to zero and the theorem easily follows.
\end{Dem}

As an easy consequence of Lemmas \ref{Lemma:1ModicaReferee} and \ref{Lemma:HoleAproximation} we have the following isovolumic approximation lemma.
\begin{Lemme}\label{Lemma:Smoothisovolumic} Let $\Omega\in\tilde{\tau}_M$ be a finite perimeter set with $V(\Omega)<+\infty$, $V(\Omega), V(\Omega^c)>0$, where $\Omega^c:=M\setminus\Omega$. Then there exists a sequence $\Omega_{k}\in\tau_{M}$ such that $V(\Omega_k)=V(\Omega)$ and $\Omega_k$ converges to $\Omega$ in the sense of finite perimeter sets. 
\end{Lemme}
\begin{Dem}[of Lemma \ref{Lemma:Smoothisovolumic}]
Let us assume that $\Omega\in\tilde{\tau}_M$ is bounded, then for any arbitrary $\varepsilon>0$, the Lemma  $\ref{Lemma:1ModicaReferee}$ applied to the finite perimeter set $\tilde{\Omega}$ given by Lemma \ref{Lemma:HoleAproximation} applied to $\Omega$, permits to find $\tilde{\Omega}_{\varepsilon}\in\tau_M$ such that $V(\tilde{\Omega}_{\varepsilon})=V(\tilde{\Omega})=V(\Omega)$ and $$V(\tilde{\Omega}_{\varepsilon}\Delta\tilde{\Omega})\leq\frac{\varepsilon}{2},$$ 
$$|\P(\tilde{\Omega}_{\varepsilon})-\P(\tilde{\Omega})|\leq\frac{\varepsilon}{2}.$$
These last two inequalities combined with \eqref{Eq:Equivalence} and \eqref{Eq:Equivalence1} imply that
\begin{equation}
V(\tilde{\Omega}_{\varepsilon}\Delta\Omega)\leq\varepsilon,
\end{equation} 
\begin{equation} 
|\P(\tilde{\Omega}_{\varepsilon})-\P(\Omega)|\leq\varepsilon.
\end{equation}
\end{Dem}

Now we are ready to prove Theorem $\ref{Thm:Equivalence}$.

\begin{Dem}[of Theorem $\ref{Thm:Equivalence}$]
Taking into account Remark $\ref{Rem:Trivialinequality}$, it is easy to check that to prove the theorem, it is enough to show the nontrivial inequality $I_M(v)\leq\tilde{I}_M(v)$ for every $v\in[0, V(M)[$.  To this aim, let us consider $\varepsilon>0$ and $\Omega\in\tilde{\tau}_M$, with $V(\Omega)=v$. By Lemma $\ref{Lemma:Smoothisovolumic}$ there is a sequence $\Omega_k\in\tau_M$ such that $V(\Omega_k)=v$, and $(\Omega_k)$ converges to $\Omega$ in the sense of finite perimeter sets. In particular we have that $\lim_{k\to+\infty}\P(\Omega_k)=\P(\Omega)$. On the other hand by definition we have that 
$I_M(v)\leq\P(\Omega_k)$ for every $k\in\mathbb{N}$. Passing to limits leads to have 
\begin{equation}\label{Eq:Thm1}
 I_M(v)\leq\P(\Omega),
\end{equation} 
for every $\Omega\in\tilde{\tau}_M$ with $V(\Omega)=v$. Taking the infimum in \eqref{Eq:Thm1} when $\Omega$ runs over $\tilde{\tau}_M$ keeping $V(\Omega)$ fixed and equal to $v$, we infer that $I_M(v)\leq\tilde{I}_M(v)$. This completes the proof.
\end{Dem}

In the remaining part of this section we prove Corollary $\ref{CorRes:Uppersemicontinuity}$.

\begin{Dem}[of Corollary $\ref{CorRes:Uppersemicontinuity}$] In view of Theorem $\ref{Thm:Equivalence}$ we actually prove that $\tilde{I}_M$ is upper semicontinuous. For any $v\in]0, V(M)[$ and any $\varepsilon>0$, consider a finite perimeter set $\Omega$ such that $V(\Omega)=v$ and $\P(\Omega)\leq\frac{\varepsilon}{4}$. We then apply Lemma \ref{Lemma:HoleAproximation} to it, which gives us $\Omega_1$ such that $V(\Omega_1)=v$, $\P(\Omega_1)\leq\tilde{I}_M(v)+\frac{\varepsilon}{2}$, and a $\bar{v}=\bar{v}_{\Omega_{1,\varepsilon}}$ such that for any $w\in]v-\bar{v}, v+\bar{v}[$ there exists $\Omega_2\in\tilde{\tau}_M$ satisfying $V(\Omega_2)=w$ and $\P(\Omega_2)\leq I_M(v)+f(|w-v|)+\frac{3\varepsilon}{4}$, where $f$ is given by \eqref{Eq:defoff}. By the very definition of isoperimetric profile we have immediately that 
$$I_M(w)\leq I_M(v)+f(|w-v|)+\frac{3\varepsilon}{4}.$$
Now, the function $f$ depends only on $\Omega_1$, satisfies $f(0)=0$ and is continuous at $0$. So there exists $v_1\in]0, \bar{v}[$ such that $f(|w-v|)\leq\frac{\varepsilon}{4}$ for every $w\in]v-v_1, v+v_1[$, which gives the upper semicontinuity in $v$. By the arbitrariness of $v$ the corollary readily follows.
\end{Dem}

\section{Local H\"older continuity of $I_M$ in bounded geometry}\label{Sec:Continuity}  


For the needs of the proof of Theorem $\ref{MainLemma}$ we restate here a version of Lemma $2.5$ of \cite{Nar12} that we will use in the sequel.
\begin{Lemme}[Lemma $2.5$ of \cite{Nar12}]\label{Lemma:NarAsian}
There is a constant $c=c(n,k)$, with $0<c<1$ such that for any Riemannian manifold $M^n$ with bounded geometry, any radius $0<r\leq 1$, any set $D\in\tilde{\tau}_M$ with $V(D)<+\infty$, there is a point $p\in M$ such that 
\begin{equation}
V(B(p,r)\cap D)\geq c\min\{v_0r^n, \left(\frac{V(D)}{\P(D)}\right)^n\}.
\end{equation} 
\end{Lemme} 

The proof of the preceding Lemma is essentially the same as in Lemma $2.5$ of \cite{Nar12}. 

Now we can start the proof of Theorem $\ref{MainLemma}$.

\begin{Dem}[of Theorem $\ref{MainLemma}$] As a preliminary remark we observe that it is enough to prove the theorem thinking to the definition of $\tilde{I}_M$ when it is more useful for our reasoning. Let $\varepsilon\in]0,1]$. By Theorem $\ref{Thm:Equivalence}$ we can get $\Omega\in\tau_M$ with $V(\Omega)=w$ and $\P(\Omega)\leq I_M(w)+\varepsilon$. When $M$ is not compact, there exists a ball $B(x_2, 1)$ not intersecting $\Omega$ (that could be chosen compact). Then for every $v'\in]w,w+v_0[$ there exists $r_{v'}\leq 1$ such that $\Omega_1=\Omega\mathring{\cup} B_M(x_2, r_{v'})$ satisfies $V(\Omega_1)=v'$ and 
\begin{small}
\begin{equation}\label{Eq:ProofMain}
I_M(v')\leq\P(\Omega_1)\leq\P(\Omega)+\P(B_M(x_2, r_{v'}))\leq I_M(w)+\varepsilon+C(n,k)r_{v'},
\end{equation}
\end{small}
where the last inequality comes from the spherical Bishop-Gromov's theorem (which asserts that when $Ric_g\geq (n-1)kg$ the area of spheres are less than the area of corresponding spheres in space form of constant curvature $k$) 
and from the value of the area of the spheres in constant curvature. Since by Bishop-Gromov's Theorem we have $\frac{v_0r_{v'}^n}{C_1(n,k)}\leq V(B(x_2, r_{v'})=v'-v$, Inequality \eqref{Eq:ProofMain} gives us
\begin{equation}\label{Eq:ProofMain0}
I_M(v')\leq I_M(v) +\varepsilon+C_2(n,k)\left(\frac{v'-v}{v_0}\right)^{\frac{n-1}{n}}.
\end{equation}
The case $v'\leq v$ needs more work. Let us apply Lemma \ref{Lemma:NarAsian} to $\Omega$, we get for any $v'\in]v-v_1, v[$, where $v_1=c\min\left\{v_0,\left(\frac{v}{I_M(v)+\varepsilon}\right)^n\right\}$, then we have 
\begin{small}
\begin{equation}
V\left(\Omega\cap B\left(p, \left(\frac{v-v'}{cv_0}\right)^{\frac{1}{n}}\right)\right)\geq\min\left\{v-v', c\left(\frac{v}{I_M(v)+\varepsilon}\right)^n\right\}=v-v',
\end{equation} 
\end{small}
and so there exists a $r_{v'}\leq\left(\frac{v-v'}{cv_0}\right)^{\frac{1}{n}}$ such that $\Omega_2:=\Omega\setminus B(p,r_{v'})$ has volume $v'$ and so, by the spherical Bishop-Gromov's Theorem, we get
\begin{small}
\begin{equation}\label{Eq:ProofMain1}
I_M(v')\leq\P(\Omega_2)\leq\P(\Omega)+\P(B_M(p,r_{v'}))\leq I_M(v)+\varepsilon+C_2(n,k)\left(\frac{v'-v}{v_0}\right)^{\frac{n-1}{n}}.
\end{equation}
\end{small}
Now, we can let $\varepsilon$ tends to $0$ in $\eqref{Eq:ProofMain0}$ and $\eqref{Eq:ProofMain1}$. If we have $v'\leq v$, then we get the result combining $\eqref{Eq:ProofMain1}$ and $\eqref{Eq:ProofMain0}$ where we exchange $v$ and $v'$. If $v\leq v$, we first control $I_M(v')$ by $I_M(v)$ using $\eqref{Eq:ProofMain0}$ and then apply $\eqref{Eq:ProofMain1}$ with $v$ and $v'$ exchanged. Combined with $\eqref{Eq:ProofMain0}$ we conclude the proof in the case $V(M)=+\infty$. If $V(M)<+\infty$ we can just take as $\Omega$ an isoperimetric region of volume $v$ (which exists always), then apply the arguments leading to \eqref{Eq:ProofMain1} to $M\setminus\Omega$ and consider as a competitor the finite perimeter set $\Omega':=\Omega\cup B_M(p,r_{v'})$, then it is straightforward to adapt the preceding arguments to conclude the proof. 
\end{Dem}

At this point, we are ready to prove Corollary $\ref{Main}$.

\begin{Dem}
Lemma $3.5$ \cite{MJ} states that whenever $(M,g)$ have $Ric_g\geq (n-1)k$ then the perimeter of a geodesic ball in $M$ enclosing volume $v$, have no more perimeter than a geodesic ball in $\mathbb{M}_k^n$ enclosing the same volume, this is used to prove Proposition $3.2$ of \cite{MonNar} which states that if $(M^n,g)$ is a complete Riemannian manifold with $Ric_g\geq (n-1)k$, then $I_M\leq I_{\mathbb{M}_k^n}$. But we know a lot about $I_{\mathbb{M}_k^n}$, for example that it is a continuous strictly increasing function and that for every $w>0$, $I_{\mathbb{M}_k^n}(w)$ is achieved by a geodesic ball enclosing volume $w$ and we will use these informations several times in the sequel. For each $v\in]0,V(M)[$ it is a trivial matter to determine $\eta_v>0$ such that $[v-\eta_v, v+\eta_v]\subset]0, V(M)[$ (for example to put $\eta_v=\min(\frac{v}{2},\frac{V(M)-v}{2})$ it is sufficient for our purposes). Put 
$$\delta:=\frac{1}{2}\min\left\{\eta_v, \frac{1}{C(n,k)}\min\left(v_0,\left(\frac{v-\eta_v}{I_{\mathbb{M}^n_k}(v+\eta_v)+C(n,k)}\right)^n\right)\right\}.$$ 
It is easy to check that $\delta=\delta(n,k,v_0, V(M),v)$.
Using Theorem $\ref{MainLemma}$ we obtain the validity of $\eqref{Eq:ResMainStatement}$ for every $v_1,v_2\in]v-\delta,v+\delta[$. To show the local uniform $\frac{n-1}{n}$-H\"older continuity away from zero we set $\delta':=\inf_{v\in[\bar{v},V(M)[}\{\delta(n,k,v_0,v)\}=\delta'(n,k,v_0)$. It is easy to see that $\delta'>0$ because $v\mapsto\delta(n,k,v_0,v)$ is a continuous function of $v$. Readily follows that $\eqref{Eq:ResMainStatement}$ holds for every $v_1,v_2\in[\bar{v}, V(M)[$ satisfying $|v_1-v_2|\leq\delta'$. Furthermore, if we assume that $V(M)<+\infty$ we can divide the interval $[\bar{v}, V(M)-\bar{v}]$ in a finite number of interval whose length is less that $\delta'$. Then it is straightforward to prove that for all $v_1,v_2\in[\bar{v}, V(M)-\bar{v}]$ we have $$|I_M(v_1)-I_M(v_2)|\leq\left(\left[\frac{V(M)-2\bar{v}}{\delta'}\right]+1\right)C(n,k)\left(\frac{|v_1-v_2|}{v_0}\right)^{\frac{n-1}{n}}.$$ To finish the proof it is enough to remark that for every $v\in[a,b]$ it holds $$\delta_v>\delta(n,k,v_0, a, b)=\frac{1}{C(n,k)}\min\left(v_0, \left(\frac{a}{I_{\mathbb{M}^n_k}(b)+C(n,k)}\right)^{n}\right)>0,$$
which ensures that $I_M$ is uniformly locally continuous on $[a,b]$. With this in mind it is a standard task to conclude the global $\frac{n-1}{n}$-H\"older continuity of $I_M$ and to complete the proof. 
\end{Dem}
      \markboth{References}{References}
      \bibliographystyle{alpha}
      \bibliography{HolderContinuityAbraham2}
      \addcontentsline{toc}{section}{\numberline{}References}
      
     \emph{Abraham Mu\~noz Flores\\ Departamento de Matem\'atica\\ Instituto de Matem\'atica\\ UFRJ-Universidade Federal do Rio de Janeiro, Brasil\\ email: abraham@im.ufrj.br\\
     and Departamento de Geometria e Representa\c{c}\~ao Gr\'afica\\
     Instituto de Matem\'atica e Estat\'istica\\
     UERJ-Universidade Estadual do Rio de Janeiro\\
     email: abraham.flores@ime.uerj.br}\\
     
      \emph{Stefano Nardulli\\ Departamento de Matem\'atica\\ Instituto de Matem\'atica\\ UFRJ-Universidade Federal do Rio de Janeiro, Brasil\\ email: nardulli@im.ufrj.br} 
\end{document}